\documentclass[a4paper,11pt]{article}
\usepackage{authblk}
\usepackage[english]{babel} 
\usepackage[centertags,intlimits]{amsmath}
\usepackage{amsfonts}
\usepackage{amsthm}
\usepackage{amssymb}
\usepackage{pgf, tikz, pgfplots, float, graphicx, tikz-cd}
\usetikzlibrary{decorations.pathreplacing,calligraphy, arrows}
\usepackage{comment}
\usepackage{cite}

\newtheorem{theorem}{Theorem}[section]
{Proposition}

{Lemma}
{Corollary}

\numberwithin{equation}{theorem}
\newcommand{\dd}[1]{\textbf{\textit{#1}}}

\begin{document}

\title{A Note on Distance-Fall Colorings}

\author[1]{Wayne Goddard}
\author[2,3]{Sonwabile Mafunda}
\affil[1]{Clemson University, USA}
\affil[2]{Soka University of America, USA}
\affil[3]{University of Johannesburg, South Africa}

\maketitle

\begin{abstract}
We say a proper coloring  of a graph is distance-$k$ fall if every vertex is within distance $k$ 
of at least one vertex of every color. We show that 
if $G$ is a connected graph of order at least $3$ that is $3$-colorable, then
it has a distance-2 fall 3-coloring. Further, for every integer $k\ge 2$,
if $T$ is a tree of order at least $k$, then $T$ has a $k$-coloring 
such that every vertex is within distance $k-1$ of every color.
This proves an old conjecture of Beineke and Henning that every tree
of order $n$ has an independent distance-$d$-dominating set of size at most $n/(d + 1)$.
\end{abstract}


MSC-class: 05C12; 05C15; 05C69

\section{Introduction}
A \dd{fall coloring} is a proper coloring of the vertices of a graph such that every vertex sees every color. 
That is, for each vertex $v$ its closed neighborhood $N[v]$ contains all colors.
Not all graphs have a fall coloring. The simplest example is the $5$-cycle.
In this paper we consider proper colorings where every vertex is ``near'' to every color.
We say a coloring is \dd{distance-$k$ fall} if every vertex is within distance $k$ 
of at least one vertex of every color. For example, 
for a graph of diameter $2$, every proper coloring is distance-$2$ fall;
and every odd cycle has a $3$-coloring that is distance-$2$ fall.

The terminology ``\dd{fall coloring}'' was introduced in the 2000 paper by Dunbar et al.~\cite{DHHJKLR}. 
But the concept is older, having been studied before as partitioning a graph into independent dominating sets;
see the references of~\cite{DHHJKLR}. Recall that a set $S$ of vertices is \dd{independent} if no two elements of $S$ are joined by an edge.
A set $S$ is \dd{dominating} if every vertex is either in $S$ or adjacent to at least one vertex of $S$.
An \dd{independent dominating} set is one that is both independent and dominating.
More generally, a set $S$ is \dd{distance-$k$ dominating} if every vertex is within distance $k$ of at least
one vertex of $S$. Thus a distance-$k$ fall coloring is a partition of the vertex set into independent distance-$k$ dominating sets.
The parameter \dd{independent distance-$k$ domination}, that is, the  
minimum size of an independent distance-$k$ dominating set,
was originally studied by Beineke and Henning~\cite{BH}.

 \section{Distance-Fall Colorings and Chromatic Number}\label{dfc&cn}

It is immediate that  if a (connected nontrivial) graph is bipartite, then the bipartite coloring is a fall coloring.
So the next question to consider is $3$-colorable graphs.
For a coloring, we say a vertex is \dd{$d$-good} if every color appears within distance $d$ of the vertex; otherwise the vertex
is \dd{$d$-bad}. The following theorem shows that if the graph is $3$-colorable then 
 there is a proper coloring with $3$ colors such that every vertex is $2$-good.

\begin{theorem} \label{t:threeGood}
If $G$ is a connected graph of order at least $3$ that is $3$-colorable, then $G$ has a distance-$2$ fall $3$-coloring.
\end{theorem}

{\bf Proof:} 
Suppose to the contrary, that $G$ does not admit a distance-$2$ fall $3$-coloring. 
Take a proper $3$-coloring of $G$ that minimizes the number of $2$-bad vertices. 
Without loss of generality, choose a $2$-bad vertex $v \in V(G)$. 
Then $N(v)$ is monochromatic.

If there exists a vertex $w_1 \in N(v)$ with $\deg(w_1)=1$, 
then by recoloring $w_1$ with the missing color in $N[v]$, 
we obtain a proper $3$-coloring of $G$ with strictly fewer $2$-bad vertices, 
contradicting the minimality of the original coloring.

Therefore, assume $\deg(w) \geq 2$ for all $w \in N(v)$. 
Since $v$ is $2$-bad, we have $N(v) \cap N(w) = \emptyset$ and $N(w)$ is monochromatic for every $w \in N(v)$. 
That is, $v$ and all vertices in its second neighborhood share the same color.

Now recolor $v$ with the missing color in $N[v]$. 
Vertex $v$ then becomes $2$-good, and every neighbor of $v$ remains $2$-good. 
Moreover, no previously $2$-good vertex becomes $2$-bad under this recoloring. 
Thus, the total number of $2$-bad vertices decreases, 
contradicting the minimality of the original coloring.

Therefore, there exists a proper $3$-coloring of $G$ with zero $2$-bad vertices.

\hfill$\Box$\\

It is perhaps interesting to note that an independent distance-$2$ dominating set is related
to an independent isolating set. An isolating (or vertex-edge dominating) set $S$ can be 
viewed as a $2$-dominating set $S$ with the added condition that if vertex $v$ is at distance $2$ 
from $S$, then all its neighbors are at distance $1$ from $S$. An independent 
isolating set is an isolating set that is also independent. See for example
\cite{BCHH,BGindep,FK}. 
It is easy to see that Theorem~\ref{t:threeGood} does not extend to a partition into three independent
isolating sets (for example $C_5$). However, in \cite{BGindep} it was shown that 
it almost does: specifically, if $G$ is $3$-colorable there exist three sets $A$, $B$, and $C$ of vertices,
such that each is an independent isolating set, their union is all vertices, 
and such that $A$ and $B$ 
overlap in at most one vertex, while $C$ is disjoint from $A \cup B$.\\

It is, however, unclear if the Theorem~\ref{t:threeGood} generalizes to larger number of colors.
The next step would be to determine whether every connected $4$-colorable
graph of order at least $4$ has a $4$-coloring such that every vertex is within distance three of all colors.
This is best possible because of path--complete graph, meaning the graph that consists of a path and a complete graph joined
by an edge.
However, in the case of trees, the above theorem does generalize:

\begin{theorem}  \label{t:treeGood}
For integer $k\ge 2$,
if $T$ is a tree of order at least $k$, then $T$ has a $k$-coloring 
such that every vertex is within distance $k-1$ of every color.
\end{theorem}

{\bf Proof:} 
Let $T$ be a tree of diameter $d$. If $d < k$, then any coloring that uses each color at least once automatically has the desired property. 

So, assume $d \geq k$. Let $P := v_1, v_2, \ldots, v_d$ be a diametral path in $T$.  
Define a coloring $c : V(P) \to \{0, 1, 2, \ldots, k-1\}$ on the vertices of $P$ such that $c(v_i) = i \pmod{k}$.

Now, let $e = uv$ be a central edge of $P$, and let $T_u$ and $T_v$ be the components of $T - e$ containing $u$ and $v$, respectively.  
If $T \setminus P$ denotes the forest induced by the vertices of $T$ not contained in $P$, then color every vertex of $T \setminus P$ such that, in each of the components $T_u$ and $T_v$, any two vertices at the same distance from $u$ or $v$ receive the same color.

We claim that this coloring has the desired property.  
If a vertex $w$ lies within distance $(k-1)/2$ of the edge $e$, then $w$ is within distance $k-1$ of each color on $P$.  
Otherwise, the unique path from $w$ to $P$ enters $P$ at some vertex and continues along $P$ via $e$; the first $k$ vertices on this path all have distinct colors.  
Hence the theorem follows. \hfill$\Box$\\

Theorem~\ref{t:threeGood} is  equivalent to saying that the vertices of a $3$-colorable graph can be partitioned into 
three (disjoint) independent distance-$2$ dominating sets. And thus the independent distance-$2$ domination
number of a $3$-colorable graph is at most $n/3$, where $n$ is the order.  
This generalizes Theorem~2 of \cite{BH}, which proved the bound for trees.
Furthermore,
Theorem~\ref{t:treeGood} shows that the independent distance-$k$ domination number of a
tree is at most $n/(k+1)$, which establishes the Conjecture at the end of~\cite{BH}.\\

We note that Theorem~\ref{t:treeGood} and thus the conjecture in~\cite{BH} was recently also proved by Bujt\'as et al.~\cite{BCKZ}. Indeed, they showed that 
Theorem~\ref{t:treeGood} generalizes to bipartite graphs.\\

There is a version of Theorem~\ref{t:threeGood} for general graphs with 
a slightly weaker distance condition. A \dd{partial coloring}
means a proper coloring where only some of the vertices are colored.

\begin{theorem} \label{t:threeGoodGeneral}
If $G$ is a connected graph of order at least $3$, then $G$ has a partial $3$-coloring
such that every vertex is within distance $3$ of every color.
\end{theorem}

{\bf Proof:} 
Consider any partial $3$-coloring that maximizes the number of colored vertices.
Note that any uncolored vertex $v$ has neighbors of every color, since it is not possible to color $v$.
Thus, out of all partial $3$-colorings with the maximum number of colored vertices, take the 
one with the minimum number of $2$-bad vertices and let $v$ be such a vertex.

We know $v$ is colored. Furthermore, all of its neighbors are colored. If some vertex $w$ at distance $2$ 
from $v$ is uncolored, then $v$ is $3$-good since $N(w)$ contains all colors. So assume every vertex within distance two of $v$ is
colored. Then, as in the proof of Theorem~\ref{t:threeGood}, one can re-color either $v$ or 
a neighbor of $v$ so that $v$ and its neighbors are $2$-good. Note that if this recoloring
makes an uncolored vertex bad, then we have a contradiction of the original requirement
that the maximum number of vertices were colored. So every vertex that is $2$-bad is
$3$-good; or in other words, every vertex is $3$-good.\hfill$\Box$\\

This strengthens the first result in \cite{BH} for the case that $k=3$.\\ 

\noindent We conclude with a brief comment about graph operations.
Kaul and Mitillos~\cite{KM} showed that if a graph $G$ has a fall $k$-coloring and a graph $H$ has a $k$-coloring,
then the cartesian product $G \Box H$ has a fall $k$-coloring. The same idea works here:
if $G$ has a distance-$d$ fall $k$-coloring $f_G$ and $H$ has a $k$-coloring $f_H$,
then the cartesian product $G \Box H$ has a distance-$d$ fall $k$-coloring.
As per usual in the cartesian product, consider both colorings as assigning
integers in the range $1$ to $k$ and take the coloring of vertex $(g,h)$ in 
the product to be the sum $f_G(g) + f_H(h)$ modulo $k$.

We note further that if $G$ has a distance-$d_G$ fall $k_G$-coloring and 
$H$ has a distance-$d_H$ fall $k_H$-coloring, then the cartesian product $G \Box H$ has 
a distance-$(d_G+d_H)$ fall $k_G k_H$-coloring. Simply take the coloring of vertex $(g,h)$ in 
the product to be the ordered pair $( f_G(g) , f_H(h) )$.\\
Similar results can be shown for other products.

\section*{Acknowledgment}
The authors would like to thank Ben Gobler for suggesting the concept.


\end{document}